\documentclass[twocolumn,11pt]{article}
\usepackage{times}

\newtheorem{definition}{Definition}
\newtheorem{remark}[definition]{Remark}

\newtheorem{theorem}[definition]{Theorem}
\newtheorem{proposition}[definition]{Proposition}
\newtheorem{corollary}[definition]{Corollary}
\newtheorem{example}[definition]{Example}

\newcommand{\eop}{\hfill $\sqcap\!\!\!\!\sqcup$} % end of proof
\begin{document}
\global\def\refname{{\normalsize \it References:}}

\baselineskip 12.5pt
\title{\LARGE \bf Special types of intuitionistic fuzzy left $\bf h$-ideals of hemirings}

\date{}

\author{\hspace*{-10pt}\begin{minipage}[t]{2.7in} \normalsize \baselineskip 12.5pt
\centerline{WIESLAW A. DUDEK}
 \centerline{Wroclaw University of
Technology}
 \centerline{Institute of Mathematics and Computer Science}
 \centerline{Wyb. Wyspianskiego 27, 50-370 Wroclaw}
 \centerline{POLAND}
  \centerline{dudek@im.pwr.wroc.pl}
  \end{minipage}
\\ \\ \hspace*{-10pt}
\begin{minipage}[b]{6.9in} \normalsize
\baselineskip 12.5pt {\it Abstract:} Characteristic, normal and
completely normal intuitionistic fuzzy left $h$-ideals of
hemirings are described.
\\ [4mm] {\it Key--Words:}
Hemiring, intuitionistic fuzzy left $h$-ideal, normal fuzzy set,
fuzzy characteristic.
\end{minipage}
\vspace{-10pt}}

\maketitle

\thispagestyle{empty} \pagestyle{empty}
\section{Introduction}
\label{S1} \vspace{-4pt} Hemirings, as semirings with zero and
commutative addition, appear in a natural manner in some
applications to the theory of automata and formal languages (see
\cite{Aho, JG, Sim}). It is a well known result that regular
languages form so-called star semirings. According to the well
known  theorem of Kleene, the languages, or sets of words,
recognized by finite-state automata are precisely those that are
obtained from letters of input alphabets by the application of the
operations: sum (union), concatenation (product), and Kleene star
(Kleene closure). If a language is represented as a formal series
with the coefficients in a Boolean hemiring, then the Kleene
theorem can be well described by the Kleen-Sch\"utzenberger
theorem. Moreover, if the coefficient hemiring is a field, then a
series is rational if and only if its syntactic algebra (see
\cite{JG, Sim, Wech} for details) has a finite rank. Many other
applications with references can be found in \cite{JG} and in a
guide to the literature on semirings and their applications
\cite{Gla}.

Ideals of hemirings play a central role in the structure theory
and are useful for many purposes. However, they do not in general
coincide with the usual ring ideals. Many results in rings
apparently have no analogues in hemirings using only ideals.
Henriksen defined in \cite{Hen} a more restricted class of ideals
in semi\-rings, which is called the class of $k$-ideals, with the
property that if the semiring $R$ is a ring then a complex in $R$
is a $k$-ideal if and only if it is a ring ideal. Another more
restricted, but very important, class of ideals, called now
$h$-ideals, has been given and investigated by Izuka \cite{Izu}
and La Torre \cite{LaT}. It interesting that the regularity of
hemirings can be characterized by fuzzy $h$-ideals \cite{ZD}.
General properties of fuzzy $k$-ideals are described in \cite{DB,
Gos, Kim}. Other important results connected with fuzzy $h$-ideals
in hemirings were obtained in \cite{Jun} and \cite{Zhan}.

\section{Preliminaries}\vspace{-4pt}

By a {\it semiring} is mean an algebraic system $(R,+,\cdot)$
consisting of a nonempty set $R$ together with two binary
operations on $R$ called addition and multiplication (denoted in
the usual manner) such that $(R,+)$ and $(R,\cdot )$ are
semigroups satisfying for all $x,y,z\in R$ the following
distributive laws

\vspace{2mm} \centerline{$x(y+z)=xy+xz$ \ \ and \ \
$(x+y)z=xz+yz$.}

By a {\it zero} we mean an element $0\in R$ such that $0x=x0=0$
and $0+x=x+0=x$ for all $x\in R$. A semiring with zero and a
commutative semigroup $(R,+)$ is called a {\it hemiring}.

A nonempty subset $A$ of $R$ is said to be a {\it left ideal} if
it is closed with respect to the addition and such that
$RA\subseteq A$. A left ideal $A$ is called a {\it left $h$-ideal}
(cf. \cite{Hen}) if for any $x,z\in R$ and $a,b\in A$ from
$x+a+z=b+z$ it follows $x\in A$.

By a {\it fuzzy set} of a hemiring $R$ we mean any mapping $\mu$
from $R$ to $[0,1]$. For any mapping $f$ from $R$ to $S$ we can
define in $R$ a new fuzzy set $\mu^f$ putting $\mu^f(x)=\mu(f(x))$
for all $x\in R$. Clearly $\mu^f(x_1)=\mu^f(x_2)$ for $x_1,x_2\in
f^{-1}(x)$.

For each fuzzy set $\mu$ in $R$ and any $\alpha\in [0,1]$ we
define two sets

\medskip
\centerline{$U(\mu,\alpha)=\{x\in R\ |\ \mu(x)\geq\alpha\}$,}

\smallskip
\centerline{$L(\mu,\alpha)=\{x\in R\ |\ \mu(x)\leq\alpha\}$,}

\medskip\noindent which are called an {\it upper} and {\it lower
level cut} of $\mu$ and can be used to the characterization of
$\mu$. The {\it complement} of $\mu$, denoted by $\overline{\mu}$,
is the fuzzy set on $R$ defined by $\overline{\mu}(x)=1-\mu(x)$.

A fuzzy set $\mu$ of a hemiring $R$ is called a {\it fuzzy left
$h$-ideal } (cf. \cite{Jun}) if for all $a,b,x,z\in R$ the
following three conditions hold:
\[\begin{array}{lll}
\mu(x+y)\ge\min\{\mu(x), \mu(y)\},\\[3pt]
\mu(xy)\ge\mu(y),\\[3pt]
x+a+z=b+z\longrightarrow \mu(x)\ge\min\{\mu(a),\mu(b)\}.
\end{array}
\]

As an important generalization of the notion of fuzzy sets,
Atanassov introduced in \cite{At1} the concept of an {\it
intuitionistic fuzzy set} (IFS for short) defined as objects
having the form:
\[
A=(\mu_A, \lambda_A)=\{(x,\mu_A(x),\lambda_A(x))\,|\,x\in R \},
 \]
where the fuzzy sets $\mu_A$ and $\lambda_A$ denote the {\it
degree of membership} (namely $\mu_A(x)$) and the {\it degree of
nonmembership} (namely $\lambda_A(x)$) of each element $x\in R$ to
the set $A$ respectively, and $0 \leq \mu_A(x) + \lambda_A(x)\leq
1$ for all $x\in R$.

According to \cite{At1}, for every two intuitionistic fuzzy sets
$A=(\mu_A, \lambda_A)$ and $B=(\mu_B,\lambda_B)$ in $R$, we
define: $A\subseteq B$ if and only if $\mu_A(x)\leq\mu_B(x)$ and
$\lambda_A(x)\geq\lambda_B(x)$ for all $x\in R$. Obviously $A=B$
means that $A\subseteq B$ and $B\subseteq A$.

\section{Intuitionistic fuzzy left $\bf h$-ideals}\vspace{-4pt}

\begin{definition}\label{D3.1}\rm
An IFS $A=(\mu_A,\lambda_A )$ on a hemiring $R$ is called an {\it
intuitionistic fuzzy left $h$-ideal } (IF left $h$-ideal for
short) if \vspace{-4pt}
\[\rule{-2mm}{0mm}\begin{array}{ll}
(1) \ \ \mu_A(x + y)\geq\min\{\mu_A(x),
 \mu_A(y)\},\\[3pt]
(2) \ \ \lambda_A(x + y)\leq\max\{\lambda_A(x),
 \lambda_A(y)\},\\[3pt]
(3) \ \ \mu_A(x y)\geq\mu_A(y),\\[3pt]
(4) \ \ \lambda_A(xy)\leq\lambda_A(y),\\[3pt]
(5)\\ x+a+z=b+z\longrightarrow\mu_A(x)\geq\min\{\mu_A(a),\mu_A(b)\},\\
(6) \\
x+a+z=b+z\longrightarrow\lambda_A(x)\leq\max\{\lambda_A(a),\lambda_A(b)\}
\end{array}
\]
hold for all $a,b,x,y,z\in R$.
\end{definition}

An IFS $A=(\mu_A,\lambda_A )$ satisfying the first four conditions
is called an {\it intuitionistic fuzzy left ideal}.

The family of all intuitionistic fuzzy left $h$-ideals of a
hemiring $R$ will be denoted by $IFI(R)$.

It is not difficult to see that $\mu_A(x)\leq\mu_A(0)$ and
$\lambda_A(0)\leq\lambda_A(x)$ for each $A\in IFI(R)$ and $x\in
R$.

\begin{example}\label{E01}\rm On a four element hemiring $(R,+,\cdot )$
defined by the following two tables:
\begin{center}
\begin{tabular}{c|cccc}

          + & 0 & 1 & 2 & 3  \\ \hline
          0 & 0 & 1 & 2 & 3  \\
          1 & 1 & 1 & 2 & 3  \\
          2 & 2 & 2 & 2 & 3  \\
          3 & 3 & 3 & 3 & 2
   \end{tabular}
  ~~~~~~~~~~
\begin{tabular}{c|cccc}
     $\cdot$& 0 & 1 & 2 & 3  \\ \hline
          0 & 0 & 0 & 0 & 0  \\
          1 & 0 & 1 & 1 & 1  \\
          2 & 0 & 1 & 1 & 1  \\
          3 & 0 & 1 & 1 & 1
\end{tabular}
  \end{center}
consider an IFS $A= (\mu_A, \lambda_A)$, where $\mu_A(0)=0.4$,
$\lambda_A(0)=0.2$ and $\mu_A(x)=0.2$, $\lambda_A(x)=0.7$ for all
$x\neq 0$. It is not difficult to verify that $A\in IFI(R)$.
\end{example}

\begin{example}\label{E02} Let $N$ be the set of all non-negative
integers an let
\[
\mu(x)=\left\{\begin{array}{cll}
 1& {\rm if } \ \ x\in\langle 4\rangle,\\[4pt]
 \frac{1}{2}& {\rm if } \ \ x\in\langle 2\rangle -\langle 4\rangle,\\[4pt]
0& {\rm otherwise},
\end{array}\right.
\]
where $\langle n\rangle$ denotes the set of all non-negative
integers divided by $n$. Then $(N,+,\cdot)$ is a hemiring and
$A=(\mu,\overline{\mu})$ is its IF left $h$-ideal.
\end{example}

 The following results can be proved by the verification of the corresponding axioms.

\begin{proposition}\label{P3.5}
A fuzzy set $\mu_A$ is a fuzzy left $h$-ideal of $R$ if and only
if $A=(\mu_A,\overline{\mu_A})$ is an IF fuzzy left $h$-ideal of
$R$.
\end{proposition}

\begin{proposition}\label{P3.6}
An IFS $A=(\mu_A,\lambda_A)$ is an IF left $h$-ideal of $R$ if and
only if $\mu_A$ and $\overline{\lambda_A}$ are fuzzy left
$h$-ideals of $R$.
\end{proposition}

\begin{proposition}{\rm \cite{Dud}}\label{P3.7}
Let $A$ be a nonempty subset of a hemiring $R$. Then an IFS
$(\mu_A,\lambda_A)$ defined by
\[
\begin{array}{lll}
\mu_A(x)=\left\{\begin{array}{ll}
 \alpha_2& {\rm if } \ \ x\in A,\\[2pt]
 \alpha_1& {\rm for} \ \ x\notin A,
 \end{array}\right.\\[16pt]
\lambda_A(x)=\left\{\begin{array}{ll}
\beta_2& {\rm if } \ \ x\in A,\\[2pt]
 \beta_1& {\rm for} \ \ x\notin A,
 \end{array}\right.
\end{array}
 \]
where $0\le\alpha_1<\alpha_2\le 1$, $0\le\beta_2<\beta_1\le 1$ and
$\alpha_i+\beta_i\le 1$ for $i=1,2$, is an IF left $h$-ideal of
$R$ if and only if $A$ is a left $h$-ideal of $R$.
\end{proposition}

\begin{definition}\label{D3.10}\rm
Let $A=(\mu_A, \lambda_A)$ be an intuitionistic fuzzy set in a
hemiring $R$ and let $\alpha, \beta \in [0,1]$ be such that
$\alpha +\beta \le 1.$ Then the set
 \[
R_A^{(\alpha,\beta)} =\{x\in R \mid \alpha\le\mu_A(x), \, \,
\lambda_A(x)\le \beta\}
    \]
is called an {\it $(\alpha,\beta)$-level subset} of
$A=(\mu_A,\lambda_A)$.

The set of all $(\alpha,\beta)\in{\rm Im}(\mu_A)\times{\rm
Im}(\lambda_A)$ such that $\alpha +\beta\le 1$ is called the {\it
image of} $A=(\mu_A,\lambda_A)$.
\end{definition}

Clearly $R_A^{(\alpha, \beta)}=U(\mu_A, \alpha)\cap
L(\lambda_A,\beta),$ where $U(\mu_A,\alpha)$ and
$L(\lambda_A,\beta)$ are upper and lower level subsets of $\mu_A$
and $\lambda_A$, respectively.

\begin{theorem}{\rm \cite{Dud}}\label{T3.11}
An $IFS$ $A=(\mu_A,\lambda_A)$ is an IF left $h$-ideal of $R$ if
and only if $R_A^{(\alpha,\beta)}$ is a left $h$-ideal of $R$ for
every $(\alpha,\beta)\in{\rm Im}(\mu_A)\times{\rm Im}(\lambda_A)$
such that $\alpha +\beta\le 1,$ i.e., if and only if all nonempty
level subsets $U(\mu_A,\alpha)$ and $L(\lambda_A,\beta)$ are left
$h$-ideals of $R$.
\end{theorem}

\begin{theorem}{\rm \cite{Dud}}\label{L3.13}
Let $A= (\mu_{A}, \lambda_{A})$ be an intuitionistic fuzzy left
$h$-ideal of a hemiring $R$ and let $x\in R$. Then
$\mu_A(x)=\alpha$, $\lambda_A(x)=\beta$ if and only if $x\in
U(\mu_A,\alpha)$, $x\notin U(\mu_A,\gamma)$ and $x\in L(\lambda_A,
\beta)$, $x\notin L(\lambda_A,\delta)$ for all $\gamma>\alpha$ and
$\delta<\beta$.
\end{theorem}

\section{Characteristic IF left $\bf h$-ideals}\vspace{-4pt}

\begin{definition}\label{D4.1}\rm
A left $h$-ideal $A$ of a hemiring $R$ is  said to be {\it
characteristic }if $f(A)=A$ for all $f\in Aut(R)$, where $Aut(R)$
is the set of all automorphisms of $R$.
 \end{definition}

\begin{definition}\label{D4.2}\rm
An IFS $A= (\mu_{A}, \lambda_{A})$ of $R$ is called an {\it
intuitionistic fuzzy characteristic} if $\mu_A^f(x)= \mu_A(x)$ and
$\lambda_A^f(x)= \lambda_A(x)$ for all $x\in R$ and $f\in Aut(R)$.
\end{definition}

\begin{theorem}\label{T4.3}
$A\in IFI(R)$ is characteristic if and only if each its nonempty
level set is a characteristic left $h$-ideal of $R$.
\end{theorem}
\noindent{\bf Proof:} An IFS $A=(\mu_A,\lambda_A)$ is an IF left
$h$-ideal if and only if all its nonempty level subsets are left
$h$-ideals (Theorem \ref{L3.13}). So, we will be prove only that
$A$ is characteristic if and only if all its level subsets are
characteristic. If $A=(\mu_A,\lambda_A)$ is characteristic,
$\alpha\in{\rm Im}(\mu_A)$, $f\in Aut(R)$, $x\in U(\mu_A,\alpha)$,
then $\mu_A^f(x)=\mu_A(f(x))=\mu_A(x)\geq\alpha$, which means that
$f(x)\in U(\mu_A,\alpha)$. Thus $f(U(\mu_A,\alpha))\subseteq
U(\mu_A,\alpha)$. Since for each $x\in U(\mu_A,\alpha)$ there
exists $y\in R$ such that $f(y)=x$ we have
$\mu_A(y)=\mu_A^f(y)=\mu_A(f(y))=\mu_A(x)\geq\alpha$. Therefore
$y\in U(\mu_A,\alpha)$. Thus $x=f(y)\in f(U(\mu_A,\alpha))$. Hence
$f(U(\mu_A,\alpha))=U(\mu_A,\alpha)$. Similarly,
$f(L(\lambda_A,\beta))= L(\lambda_A,\beta)$. This proves that
$U(\mu_A,\alpha)$ and $L(\lambda_A,\beta)$ are characteristic.

Conversely, if all levels of $A=(\mu_A,\lambda_A)$ are
characteristic left $h$-ideals of $R$, then for $x\in R$, $f\in
Aut(R)$ and $\mu_A(x)=\alpha$, $\lambda_A(x)=\beta$, by Lemma
\ref{L3.13}, we have $x\in U(\mu_A,\alpha)$, $x\notin
U(\mu_A,\gamma)$ and $x\in L(\lambda_A,\beta)$, $x\notin
L(\lambda_A,\delta)$ for all $\gamma>\alpha$, $\delta<\beta$. Thus
$f(x)\in f(U(\mu_A,\alpha))=U(\mu_A,\alpha)$ and $f(x)\in
f(L(\lambda_A,\beta))= L(\lambda_A,\beta)$, i.e., $\mu_A(f(x))\geq
\alpha$ and $\lambda_A(f(x))\leq\beta$. For
$\mu_A(f(x))=\gamma>\alpha$, $\lambda_A(f(x))=\delta<\beta$ we
have $f(x)\in U(\mu_A,\gamma)=f(U(\mu_A,\gamma))$, $f(x)\in
L(\lambda_A,\delta)= f(L(\lambda_A,\delta))$, which implies $x\in
U(\mu_A,\gamma)$, $x\in L(\mu_A,\delta)$. This is a contradiction.
Thus $\mu_A(f(x))= \mu_A(x)$ and $\lambda_A(f(x))=\lambda_A(x)$.
So, $A=(\mu_A,\lambda_A)$ is characteristic.
 \eop

\begin{proposition}\label{P4.4}
Let $f:R\to S$ be a homomorphism of hemirings. If
$A=(\mu_A,\lambda_A)$ is an IF left $h$-ideal of $S$, then
$A^f=(\mu^f_A ,\lambda^f_A)$ is an IF left $h$-ideal of $R$.
\end{proposition}
\noindent{\bf Proof:} Let $x,y\in R$. Then
\[\arraycolsep=.5mm\begin{array}{rll}
\mu^f_A(x+y)&=\mu_A(f(x+y))=\mu_A(f(x)+f(y))\\[3pt]
&\geq\min\{\mu_A(f(x)),\mu_A(f(y))\}\\[3pt]
&=\min\{\mu^f_A(x),\mu^f_A(y)\},\\[8pt]
\lambda^f_A(x+y)&=\lambda_A(f(x+y))=\lambda_A(f(x)+f(y))\\[3pt]
&\leq\max\{\lambda_A(f(x)),\lambda_A(f(y))\}\\[3pt]
&=\max\{\lambda_A^f(x),\lambda^f_A(y)\},\\[8pt]
\mu^f_A(xy)&=\mu_A(f(xy))=\mu_A(f(x)f(y))\\[3pt]
&\geq\mu_A(f(y))=\alpha^f_A(y),\\[8pt]
\lambda^f_A(xy)&=\lambda_A(f(xy))=\lambda_A(f(x)f(y))\\[3pt]
&\leq\lambda_A(f(y))=\lambda^f_A(y).
\end{array}\]
If $x+a+z=b+z$, then $f(x)+f(a)+f(z)=f(b)+f(z)$, whence
\[\arraycolsep=.5mm\begin{array}{rll}
\mu^f_A(x)&=\mu_A(f(x))\geq\min\{\mu_A(f(a))\mu_A(f(b))\}\\[3pt]
&=\min\{\mu^f_A(a),\mu^f_A(b)\},\\[6pt]
\lambda^f_A(x)&=\lambda_A(f(x))\leq\max\{\lambda_A(f(a)),\lambda_A(f(b))\}\\[3pt]
&=\max\{\lambda^f_A(a), \lambda^f_A(b)\}.
\end{array}\]
This proves that $A^f=(\mu^f_A,\lambda^f_A)$ is an IF left
$h$-ideal of $R$. \eop

\begin{proposition}\label{P4.5}
Let $f:R\to S$ be an epimorphism of hemirings. If
$A^f=(\mu^f_A,\lambda^f_A)$ is an IF left $h$-ideal of $R$, then
$A=(\mu_A,\lambda_A)$ is an IF left $h$-ideal of $S$.
\end{proposition}
\noindent{\bf Proof:} Since $f$ is a surjective mapping, for
$x,y\in S$ there are $x_1,y_1\in R$ such that $x=f(x_1)$,
$y=f(y_1)$. Thus
\[\arraycolsep=.5mm
\begin{array}{rlll}
\mu_A(x+y)&=\mu_A(f(x_1)+f(y_1))=\mu_A(f(x_1+y_1))\\[3pt]
&=\mu^f_A(x_1+y_1)\geq\min\{\mu^f_A(x_1),\mu^f_A(y_1)\}\\[3pt]
&=\min\{\mu_A(x),\mu_A(y)\},
\end{array}
\]
proves that $\mu_A$ satisfies the first condition on Definition
\ref{D3.1}. In a similar way we can verify others conditions. \eop

\medskip
As a consequence of the above two propositions we obtain the
following theorem.

\begin{theorem}\label{T4.6}
Let $f:R\to S$ be an epimorphism of hemirings. Then
$A^f=(\mu^f_A,\lambda^f_A)$ is an IF left $h$-ideal of $R$ if and
only if $A=(\mu_A,\lambda_A)$ is an IF left $h$-ideal of $S$.
\end{theorem}

\section{Normal IF left $\bf h$-ideals}

\begin{definition}\label{D4.7}\rm
An IF left $h$-ideal $A=(\mu_A,\lambda_A)$ of a hemiring $R$ is
said to be {\it normal} if $A(0)=(1,0)$, i.e., $\mu_A(0)=1$ and
$\lambda_A(0)= 0$.
\end{definition}

It is clear that any IF left $h$-ideal containing a normal IF left
$h$-ideal is normal too.

\begin{theorem}\label{T4.8}
Let $A=(\mu_A, \lambda_A)\in IFI(R)$ and let $\mu^+_A(x)=
\mu_A(x)+1-\mu_A(0)$, $\lambda^+_A(x)= \lambda_A(x)-\lambda_A(0)$.
If $\mu^+_A(x)+\lambda^+_A(x)\leq 1$ for all $x\in R$, then
$A^+=(\mu^+_A,\lambda^+_A)$ is a normal IF left $h$-ideal of $R$
containing $A$.
\end{theorem}
\noindent{\bf Proof:} At first observe that $\mu^+_A(0)=1$,
$\lambda^+_A(0)=0$ and $\mu^+(x),\lambda^+(x)\in [0,1]$ for every
$x\in R$. So, $A^+=(\mu^+_A,\lambda^+_A)$ is a normal IFS.

To prove that it is an IF left $h$-ideal let $x,y\in R$. Then
\[\arraycolsep=.5mm\begin{array}{lll}
\mu^+_A(x+y) = \mu_A(x+y)+1-\mu_A(0) \\
\geq\min\{\mu_A(x),\mu_A(y)\}+1-\mu_A(0) \\
=\min\{\mu_A(x)+1-\mu_A(0),\mu_A(y)+1-\mu_A(0)\} \\
=\min\{\mu^+_A(x),\mu^+_A(y)\},\\[6pt]
\lambda^+_A(x+y)= \lambda_A(x+y)-\lambda_A(0) \\
\leq\max\{\lambda_A(x),\lambda_A(y)\}-\lambda_A(0) \\
=\max\{\lambda_A(x)-\lambda_A(0),\lambda_A(y)-\lambda_A(0)\} \\
=\max\{\lambda^+_A(x),\lambda^+_A(y)\},
\end{array}
\]
and
\[\arraycolsep=.5mm
\begin{array}{rll}
 \mu^+_A(xy)&= \mu_A(xy) +1-\mu_A(0)\\
 &\geq\mu_A(y)+1-\mu_A(0)=\mu^+_A(y),\\[6pt]
\lambda^+_A(xy)&=\lambda_A(xy)-\lambda_A(0)\\
&\leq\lambda_A(y)-\lambda_A(0)=\lambda^+_A(y).
\end{array}\]
This shows that $A^+$ is an IF left ideal of $R$. Moreover, if
$x+a+z=b+z$, then
\[\begin{array}{ll}
\mu^+_A(x) = \mu_A(x) +1-\mu_A(0) \\
\geq\min\{\mu_A(a), \mu_A(b)\} +1-\mu_A(0) \\
=\min\{\mu_A(a)+1-\mu_A(0), \mu_A(b)+1-\mu_A(0)\}  \\
=\min\{\mu^+_A(a), \mu^+_A(b)\}.
\end{array}\]
Similarly
\[\begin{array}{ll}
\lambda^+_A(x)=\lambda_A(x)-\lambda_A(0) \\
\leq\max\{\lambda_A(a), \lambda_A(b)\}-\lambda_A(0) \\
=\max\{\lambda_A(a)-\lambda_A(0), \lambda_A(b)-\lambda_A(0)\}  \\
=\max\{\lambda^+_A(a),\lambda^+_A(b)\}.
\end{array}
\]
So, $A^+=(\mu^+_A,\lambda^+_A)$ is a normal IF left $h$-ideal of
$R$. Clearly $A\subseteq A^+$. \eop

\begin{remark}\rm
In the above theorem the assumption $\mu^+_A(x)+\lambda^+_A(x)\leq
1$ is essential. Indeed, for an IF left $h$-ideal $A$ defined in
Example \ref{E01} for every $x\in R$ we have
$\mu_A^+(x)=\mu_A(x)+0.6\in [0,1]$ and
$\lambda^+_A(x)=\lambda_A(x)-0.2\in [0,1]$, but
$\mu^+_A(x)+\lambda^+_A(x)>1$ for all $x>0$. So,
$A^+=(\mu^+_A,\lambda^+_A)$ is not an IFS.
\end{remark}

\begin{corollary}\label{C4.9}
$(A^+)^+=A^+$ for any $A\in IFI(R)$. If $A$ is normal, then
$A^+=A$.
\end{corollary}

Denote by $NIFI(R)$ the set of all normal IF left $h$-ideals of
$R$. Note that $NIFI(R)$ is a poset under the set inclusion.

\begin{theorem}\label{T4.10}
A non-constant maximal element of $(NIFI(R),\subseteq)$ takes only
the values $(0,1)$ and $(1,0)$.
\end{theorem}
\noindent{\bf Proof:} Let $A=(\mu_A,\lambda_A)\in NIFI(R)$ be a
non-constant maximal element of $(NIFI(R), \subseteq)$. Then
$\mu_A(0)=1$ and $\lambda_A(0)=0$. Let $x\in R$ be such that
$\mu(x)\ne 1$. We claim that $\mu_A(x)=0$. If not, then there
exists $c\in R$ such that $0<\mu_A(c)<1$. Let $A_c=(\nu_A,\rho_A)$
be an IFS in $R$ defined by
\[\begin{array}{ll}
\nu_A(x)=\frac{1}{2}\{\mu_A(x)+\mu_A(c)\},\\[4pt]
\rho_A(x)=\frac{1}{2}\{\lambda_A(x)+\lambda_A(c)\}
\end{array}
\]
for all $x\in R$. Then clearly an IFS $A_c$ is well-defined and
\[\arraycolsep=.5mm
\begin{array}{rl}
\nu_A(0)&=\frac{1}{2}\{\mu_A(0)+\mu_A(c)\}\\[3pt]
&\geq\frac{1}{2}\{\mu_A(x)+\mu_A(c)\}=\nu_A(x),\\[6pt]
\rho_A(0)&=\frac{1}{2}\{\lambda_A(0)+\lambda_A(c)\}\\[3pt]
&\leq \frac{1}{2}\{\lambda_A(x)+\lambda_A(c)\}=\rho_A(x)
\end{array}
\]
for all $x\in R.$ For all $x,y\in R$ we have also
\[\arraycolsep=.5mm
\begin{array}{ll}
\nu_A(x+y)=\frac{1}{2}\{\mu_A(x+y)+\mu_A(c)\}\\[4pt]
\geq \frac{1}{2}\{\min\{\mu_A(x),\mu_A(y)\}+\mu_A(c)\}\\[4pt]
=\min\{\frac{1}{2}\{\mu_A(x)+\mu_A(c)\},
\frac{1}{2}\{\mu_A(y)+\mu_A(c)\}\}\\[4pt]
=\min\{\nu_A(x),\nu_A(y)\}
\end{array}
\]
and
\[\arraycolsep=.5mm
\begin{array}{rll}
\nu_A(xy)&=\frac{1}{2}\{\mu_A(xy)+\mu_A(c)\}\\[4pt]
&\geq\frac{1}{2}\{\mu_A(y)+\mu_A(c)\} =\nu_A(y).
\end{array}
\]
Analogously $\rho_A(x+y)\leq\max\{\rho_A(x),\rho_A(y)\}$ and
$\rho_A(xy)\leq\rho_A(y)$.

Moreover, if $x+a+z=b+z$, then
 \[
\begin{array}{l}
\nu_A(x)=\frac{1}{2}\{\mu_A(x)+\mu_A(c)\}\\[4pt]
\geq \frac{1}{2}\{\min\{\mu_A(a),\mu_A(b)\}+\mu_A(c)\}\\[4pt]
=\min\{\frac{1}{2}\{\mu_A(a)+\mu_A(c)\},
\frac{1}{2}\{\mu_A(b)+\mu_A(c)\}\}\\[4pt]
=\min\{\nu_A(a),\nu_A(b)\}
\end{array}
\]
and, by analogy, $\rho_A(x)\leq\max\{\rho_A(a),\rho_A(b)\}$. This
proves that $A_c\in IFI(R)$.

According to Theorem \ref{T4.8}, $A^+_c=(\nu_A^+,\rho_A^+)$, where
$\nu^+_A(x)=\nu_A(x)+1-\nu_A(0)=\frac{1}{2}\{1+\mu_A(x)\}$ and
$\rho^+_A(x)=\rho_A(x)-\rho_A(0)=\frac{1}{2}\lambda_A(x)$, belongs
to $NIFI(R)$. Clearly $A\subseteq A^+_c$.

Since $\nu^+_A(x)=\frac{1}{2}(1+\mu_A(x))>\mu_A(x)$, $A$ is a
proper subset of $A^+_c$. Obviously $\nu^+_A(a)<1=\nu^+_A(0)$.
Hence $A^+_c$ is non-constant, and $A$ is not a  maximal element
of $NIFI(R)$. This is a contradiction. Therefore $\mu_A$ takes
only two values: $0$ and $1$.

Analogously we can prove that $\lambda_A$ also takes the values
$0$ and $1$. This means that for $A$ the possible values are
$(0,0)$, $(0,1)$ and $(1,0)$. If $A$ takes these three values,
then
\[\arraycolsep=.5mm
\begin{array}{ll}
R^{(0,0)}&=\{x\in R\,|\,\mu_A(x)\geq 0, \ \ \lambda_A(x)\leq
0\}\\[3pt]
&=\{x\in R\,|\,\lambda_A(x)=0\},\\[4pt]
R^{(1,0)}&=\{x\in R\,|\,\mu_A(x)\geq 1, \ \ \lambda_A(x)\leq
0\}\\[3pt]
&=\{x\in R\,|\,\mu_A(x)=1, \ \ \lambda_A(x)=0\},\\[4pt]
R^{(0,1)}&=\{x\in R\,|\,\mu_A(x)\geq 0, \ \ \lambda_A(x)\leq 1\}=R
\end{array}
\]
are nonempty left $h$-ideals (Theorem \ref{T3.11}) such that
$R^{(1,0)}\subset R^{(0,0)}\subset R^{(0,1)}=R$. Then, according
to Proposition \ref{P3.7}, an IFS $B=(\mu_B,\lambda_B)$ defined by
\[
\begin{array}{lll}
\mu_B(x)=\left\{\begin{array}{llll}
 1& {\rm if }& \ x\in R^{(0,0)}\\[2pt]
 0& {\rm if }& \ x\notin R^{(0,0)}
 \end{array}\right.\\[16pt]
\lambda_B(x)=\left\{\begin{array}{llll}
0& {\rm if } & \ x\in R^{(0,0)}\\[2pt]
1& {\rm if }& \ x\notin R^{(0,0)}
 \end{array}\right.
\end{array}
 \]
is an intuitionistic fuzzy left $h$-ideals of $R$. It is normal.
Moreover, $\lambda_A(x)\ne 0$ for $x\in R\backslash R^{(0,0)}$.
Thus $\lambda_A(x)=1$, consequently $\mu_A(x)=0$. This implies
$A(x)=B(x)$ for $x\in R\backslash R^{(0,0)}$. For $x\in R^{(0,0)}$
we have $\lambda_A(x)=0=\lambda_B(x)$ and $\mu_A(x)\leq
1=\mu_B(x)$. Hence $A\subset B$. Since $\mu_A(x)=0<\mu_B(x)$ for
$x\in R^{(0,0)}\backslash R^{(1,0)}$, an IF left $h$-ideal $A$ is
a proper subset of $B$. This is a contradiction. So, a
non-constant maximal element of $(NIFI(R), \subseteq)$ takes only
two values: $(0,1)$ and $(1,0)$. \eop

\begin{definition}\label{D4.11}\rm
A non-constant $A\in IFI(R)$ is called {\it maximal } if $A^+$ is
a maximal element of the poset $(NIFI(R), \subseteq)$.
\end{definition}

\begin{theorem} \label{T4.12}
A maximal $A\in IFI(R)$ is normal and takes only two values:
$(0,1)$ and $(1,0)$.
\end{theorem}
\noindent{\bf Proof:} Let $A\in IFI(R)$ be maximal. Then $A^+$ is
a non-constant maximal element of  $(NIFI(R), \subseteq)$ and, by
Theorem \ref{T4.10}, the possible values of $A^+$ are $(0,1)$ and
$(1,0)$, i.e, $\mu^+_A$ takes only two values $0$ and $1$. Clearly
$\mu^+_A(x)=1$ if and only if $\mu_A(x)=\mu_A(0)$, and
$\mu^+_A(x)=0$ if and only if $\mu_A(x)=\mu_A(0)-1$. But
$A\subseteq A^+$ (Theorem \ref{T4.8}), so,
$\mu_A(x)\leq\mu^+_A(x)$ for all $x\in R$. Thus $\mu^+_A(x)=0$
implies $\mu_A(x)=0$, whence $\mu_A(0)=1$. This proves that $A$ is
normal. \eop

\begin{theorem} \label{T4.13}
A $(1,0)$-level subset of a maximal IF left $h$-ideal of $R$ is a
maximal left $h$-ideal of $R$.
\end{theorem}
\noindent{\bf Proof:} Let $S$ be a $(1,0)$-level subset of a
maximal $A\in IFI(R)$, i.e.,
\[
S=R^{(1,0)}=\{x\in R\,|\,\mu_A(x)=1\}.
 \]
It is not difficult to verify that $S$ is a left $h$-ideal of $R$.
$S\ne R$ because $\mu_A$ takes two values.

Let $M$ be a left $h$-ideal of $R$ containing $S$. Then
$\mu_S\subseteq\mu_M$. Since $\mu_A=\mu_S$ and $\mu_A$ takes only
two values, $\mu_M$ also takes these two values. But, by the
assumption, $A\in IFI(R)$ is maximal, so $\mu_S=\mu_A=\mu_M$ or
$\mu_M(x)=1$ for all $x\in R$. In the last case $S=R$, which is
impossible. So, $\mu_A=\mu_S=\mu_M$, which implies $S=M$. This
means that $S$ is a maximal left $h$-ideal of $R$. \eop

\begin{definition}\label{D4.14}\rm A normal $A\in IFI(R)$ is called
{\it completely normal} if there exists $x_0\in R$ such that
$A(x_0)=(0,1)$.
\end{definition}

Denote by ${\cal C}(R)$ the set of all completely normal $A\in
IFI(R)$. Clearly ${\cal C}(R)\subseteq NIFI(R)$.

\begin{theorem}\label{T4.15} A non-constant maximal element of
$(NIFI(R), \subseteq)$ is also a maximal element of $({\cal C}(R),
\subseteq)$.
 \end{theorem}
\noindent{\bf Proof:} Let $A$ be a non-constant maximal element of
$(NIFI(R), \subseteq)$. By Theorem \ref{T4.10}, $A$ takes only the
values $(0,1)$ and $(1,0)$, so $A(0)=(1,0)$ and $A(x_0)=(0,1)$ for
some $x_0\in R$. Hence $A\in {\cal C}(R)$. Assume that there
exists $B\in {\cal C}(R)$ such that $A\subseteq B$. It follows
that $A\subseteq B$ in $NIFI(R)$. Since $A$ is maximal in
$(NIFI(R), \subseteq )$ and since $B$ is non-constant, therefore
$A=B$. Thus $A$ is maximal element of $({\cal C}(R), \subseteq)$,
ending the proof.
 \eop

\begin{theorem}\label{T4.16} Every maximal $A\in IFI(R)$ is
completely normal.
\end{theorem}
\noindent{\bf Proof:} Let $A\in IFI(R)$ be maximal. Then by
Theorem \ref{T4.12}, it is normal and $A=A^+$ takes only two
values $(0,1)$ and $(1,1)$. Since $A$ is non-constant, it follows
that $A(0)=(1,0)$ and $A(x_0)=(0,1)$ for some $x_0\in R$. Hence
$A$ is completely normal, ending the proof.
 \eop

\medskip
Below we present the method of construction a new normal
intuitionistic fuzzy left $h$-ideal from old.

\begin{theorem}\label{T4.17}
Let $f:[0,1]\to [0,1]$ be an increasing function and let
$A=(\mu_A,\lambda_A)$ be an IFS on a hemi\-ring $R$. Then $A_f
=(\mu_{A_f},\lambda_{A_f})$, where $\mu_{A_f}(x)= f(\mu_A(x))$ and
$\lambda_{A_f}(x)= f(\lambda_A(x))$, is an IF left $h$-ideal if
and only if $A=(\mu_A,\lambda_A)$ is an IF left $h$-ideal.
Moreover, if $f(\mu_A(0))=1$ and $f(\lambda_A(0))=0$, then $A_f$
is normal.
\end{theorem}
\noindent{\bf Proof:} We will verify only the condition $(1)$. Let
$A_f=(\mu_{A_f},\lambda_{A_f})\in IFI(R)$. Then
\[
\begin{array}{rll}
f(\mu_A(x+y))&=\mu_{A_f}(x+y)\\[4pt]
&\geq\min\{\mu_{A_f}(x),\mu_{A_f}(y)\}\\[4pt]
&=\min\{f(\mu_A(x)),f(\mu_A(y))\}\\[4pt]
&=f(\min\{\mu_A(x),\mu_A(y)\}),
\end{array}
\]
i.e.,

\centerline{$f(\mu_A(x+y))\geq f(\min\{\mu_A(x),\mu_A(y)\})$}

\medskip\noindent for all $x,y\in R$. Since $f$ is increasing, must be
 \[
\mu_A(x+y)\geq \min\{\mu_A(x),\mu_A(y)\}.
 \]

Conversely, if $A=(\mu_A,\lambda_A)\in IFI(R)$, then for all
$x,y\in R$ we have
 \[\arraycolsep=.5mm
 \begin{array}{rlll}
\mu_{A_f}(x+y)&=f(\mu_A(x+y))\\[4pt]
&\geq f(\min\{\mu_A(x),\mu_A(y)\})\\[4pt]
&= \min\{f(\mu_A(x)),f(\mu_A(y))\}\\[4pt]
& = \min\{\mu_{A_f}(x),\mu_{A_f}(y)\},
\end{array}
\]
i.e., $\;\mu_{A_f}(x+y)\geq\min\{\mu_{A_f}(x),\mu_{A_f}(y)\}.$

This proves that $A_f=(\mu_{A_f},\lambda_{A_f})$ satisfies (1) if
and only if it is satisfying by $A=(\mu_{A},\lambda_{A})$.

In the same manner we can proved the analogous connections for the
axioms $(2)-(6)$. \eop

\section{Conclusion}
\label{S3} \vspace{-4pt}

In the present paper we present the basic results on IF left
$h$-ideals of hemirings. In our opinion the future study of
different types of IF ideals in hemirings and near rings can be
connected with $(1)$ investigating semiprime and prime IF
$h$-ideals; $(2)$ finding intuitionistic and/or interval valued
fuzzy sets and triangular norms. The obtained results can be used
to solve some social networks problems, automata theory and formal
languages.

\end{document}